\documentclass[
12pt,twoside, a4paper,
]{amsart}

\title[Unramified Brauer group]
{ On the unramified Brauer group\\
 of a homogeneous space
}
\author{Mikhail Borovoi}
\address{Raymond and Beverly Sackler School of Mathematical Sciences,
Tel Aviv University, 69978 Tel Aviv, Israel}
\email{borovoi@post.tau.ac.il}
\thanks{Partially supported by the
Hermann Minkowski Center for Geometry}


\usepackage[arrow,matrix]{xy}

\usepackage{mathptmx}
\usepackage{amscd}
\usepackage{amssymb}
\usepackage{eucal}
\usepackage{amsxtra}

\usepackage{verbatim}
\usepackage{url}

\DeclareSymbolFont{rsfs}{U}{rsfs}{m}{n}
\DeclareSymbolFontAlphabet{\mathcal}{rsfs}

\DeclareFontEncoding{OT2}{}{} 
\DeclareTextFontCommand{\textcyr}{\fontencoding{OT2}
    \fontfamily{wncyr}\fontseries{m}\fontshape{n}\selectfont}

\usepackage{verbatim}

\theoremstyle{plain}

\newtheorem{theorem}{Theorem}

\newtheorem{conditional-result}[theorem]{Conditional Result}

\newtheorem{theorem?}{Theorem(?)}[section]
\newtheorem{proposition?}[theorem]{Proposition(?)}
\newtheorem{lemma?}[theorem]{Lemma(?)}
\newtheorem{corollary?}[theorem]{Corollary(?)}

\newtheorem*{theorem*}{Theorem}
\newtheorem*{proposition*}{Proposition}
\newtheorem*{lemma*}{Lemma}
\newtheorem*{corollary*}{Corollary}
\newtheorem*{question*}{Question}
\newtheorem*{conjecture*}{Conjecture}
\newtheorem*{claim*}{Claim}

\theoremstyle{definition}

\newtheorem*{definition*}{Definition}
\newtheorem*{example*}{Example}

\theoremstyle{remark}

\newtheorem*{remark*}{Remark}
\newtheorem*{Remarks*}{Remarks}
\newenvironment{remarks*}{\begin{Remarks*}\nopagebreak[4]
\rule{1em}{0ex}\par 
\begin{theoremlist}}%
{\end{theoremlist}\end{Remarks*}}

\newcommand{\into}{\hookrightarrow}

\newcommand{\id}{\operatorname{id}}

\newcommand{\Cc}{{\mathbb{C}}}

\newcommand{\uu}{\mathrm{u}}

\def\red{\mathrm{red}}
\def\tor{{\mathrm{tor}}}

\def\sss{{\mathrm{ss}}}

\def\ssu{{\mathrm{ssu}}}

\DeclareMathOperator{\Br}{Br}

\def\Brnr{{{\rm Br}_{\rm nr}}\,}
\def\Brnrr{{{\rm Br}_{\rm nr}}}

\def\Br{{\rm Br\ }}
\def\Brr{{\rm Br\,}}
\def\nr{{\rm nr}}
\def\Pic{{\rm Pic\ }}

\begin{document}

\maketitle

\begin{abstract}
We give a new proof of the theorem stating that
for any connected linear algebraic group $G$
over an algebraically closed field $k$ of characteristic 0
and for any connected closed subgroup $H$ of $G$,
the unramified Brauer group of $G/H$ vanishes.
\end{abstract}

\section{Introduction}

In this note $k$ always denotes an algebraically closed field of characteristic 0.
For an irreducible algebraic variety $X$ over $k$,
we denote by $k(X)$ the field of rational functions on $X$.
We denote by $\Brnr\, k(X)$, or just by $\Brnr X$,
the unramified Brauer group of $k(X)$ with respect to $k$,
see \cite[Def.~5.3]{Mumbai}.
\smallskip

We give a new proof of the following theorem:

\begin{theorem}[{\cite[Thm.~5.1]{BDH}}]\label{thm:main}
Let $G$ be a connected linear algebraic group over an algebraically closed field $k$ of characteristic 0,
and let $H\subset G$ be a  connected closed subgroup.
Then $\Brnr\, k(G/H)=0$.
\end{theorem}

In the case when $G$ is simply connected this  is a classical result of Bogomolov \cite[Thm.~2.4]{Bog89},
see Colliot-Th\'el\`ene and Sansuc \cite[Thm.~9.13]{Mumbai}.
Bogomolov proved his theorem by a topological method in the case $k=\Cc$,
but the general case of an algebraically closed field $k$ of characteristic 0
reduces to the case $k=\Cc$ by the Lefschetz principle, see \cite{Mumbai}, beginning of \S\,9.
Our Theorem \ref{thm:main} answers affirmatively a question of Colliot-Th\'el\`ene and Sansuc
in \cite[Rem.~ 9.14]{Mumbai}
and a question after Theorem 1.4 in  the paper \cite{CTK2} by Colliot-Th\'el\`ene and Kunyavski\u\i.
Theorem \ref{thm:main} was recently proved   in the preprint \cite{BDH} of Demarche, Harari and the author
by a number-theoretical method.
Here we deduce this theorem from Bogomolov's theorem.
\medskip

\noindent
{\em Acknowledgements.}
The author is very grateful to Cyril Demarche for his helpful suggestions,
which permitted to considerably shorten and simplify the proof.
The author thanks Jean-Louis Colliot-Th\'el\`ene for useful remarks.

\section{Notation and preliminaries}\label{sec:prel}

By $k$ we always denote an algebraically closed field of characteristic 0.
Let $G$ be a connected linear algebraic group over $k$.
We use the following notation:

$G^\uu$ is the unipotent radical of $G$;

$G^\red=G/G^\uu$, it is reductive;

$G^\sss=[G^\red, G^\red]$, it is semisimple;

$G^\tor=G^\red/G^\sss$, it is a torus;

$G^\ssu=\ker[G\to G^\tor]$, it is an extension of a semisimple group $G^\sss$ by a unipotent group $G^\uu$.

Note that $G^\tor$ is the largest quotient torus of $G$.
Note also that $\Pic G=0$ if and only if $G^\ssu$ is simply connected, cf. \cite{Sansuc}, Lemma 6.9 and Remark 6.13.
\smallskip

Let $X$ be a smooth integral variety over $k$.
If $V$ is a smooth compactification of $X$ (existing by Hironaka's theorem),
then we can identify $\Brnr X$ with $\Br V$, cf.~\cite[Thm.~5.11(iii)]{Mumbai}.
We regard  $\Brnr X=\Br V$ as a subgroup of $\Br X$, cf.~\cite[Thm.~5.11(ii)]{Mumbai}.
If $f\colon X_1\to X_2$ is a morphism of smooth integral varieties defined over $k$,
one can extend $f$ to a morphism of suitable smooth compactifications
$f'\colon V_1\to V_2$, where $V_i$ is a smooth compactification of $X_i$ $(i=1,2)$,
see \cite[\S\,1.2.2]{BK} (again, one uses Hironaka's theorem).
It follows that $f$ induces a homomorphism of the unramified Brauer groups
$f^{\nr}\colon \Brnr X_2\to\Brnr X_1$ fitting into a commutative diagram
\begin{equation}\label{eq:fuctoriality}
\xymatrix{
\Brnr X_2\ar@{^{(}->}[d]\ar[r]^{f^{\,\nr}}      &\Brnr X_1\ar@{^{(}->}[d] \\
 \Br X_2\ar[r]^{f^*}                            &\Br X_1.
}
\end{equation}

\section{Reduction to a special case}\label{particular-case}

Let $G$ be a connected linear algebraic group defined over $k$,
and let $H\subset G$ be a connected closed subgroup.

\subsection{Reduction to the case $\Pic G=0$}
It is well known (see e.g. \cite[Lemma 5.2]{Bor-Crelle}) that there exists
a connected linear algebraic group $G'$ over $k$ with $\Pic G'=0$ and a connected closed subgroup $H'\subset G'$,
such  that the varieties $G/H$ and $G'/H'$ are isomorphic.
Therefore, we may and shall assume that $G$ in Theorem \ref{thm:main} satisfies $\Pic G=0$.

\subsection{Reduction to the case  $H=H^\ssu$}
Consider the subgroup $H^\ssu\subset H$.
The map $G/H^\ssu\to G/H$ is a right $H^\tor$-torsor.
Since $H^\tor$ is a split torus,
by Hilbert's Theorem 90 this torsor admits a local section.
Thus  the homogeneous space  $G/H^\ssu$ is birationally equivalent to $G/H\times_k H^\tor$,
and by \cite[Prop.~5.7]{Mumbai} we have $\Brnrr(G/H^\ssu)\cong\Brnrr(G/H)$.
Therefore, we may and shall assume in Theorem \ref{thm:main}  that $H=H^\ssu$,
i.e. $H$ is character-free.

\section{Deduction of Theorem \ref{thm:main} from Bogomolov's theorem}

Consider the map $G\to G/H$.
Since  the variety of $G$ is rational, we have $\Brnr G=0$,
and we see from diagram \eqref{eq:fuctoriality}
that
\begin{equation}\label{eq:rational-group}
\Brnrr(G/H)\subset\ker[\Brr(G/H)\to\Br G].
\end{equation}

Taking in account the results of Section \ref{particular-case},
we now assume that  $\Pic G=0$
and that $H\subset G$ is connected and character-free.
Set $G_1=G^\ssu$, then $G_1$ is simply connected because $\Pic G=0$, see \S\,\ref{sec:prel}.
Since $H$ is character-free, we have $H\subset G_1$.

Let $i\colon G_1\into G$ denote the inclusion homomorphism.
Consider the following commutative diagram of morphisms of varieties:
\begin{equation}\label{eq:diagr-G-H}
\xymatrix{
G_1\ar[d]\ar[r]^i      &G\ar[d]   \\
G_1/H\ar[r]^{i_*}      &G/H.
}
\end{equation}
By functoriality (see \S\,\ref{sec:prel}) this diagram defines a homomorphism $i^{\,\nr}\colon \Brnrr(G/H)\to\Brnrr(G_1/H)$
fitting into a commutative diagram
\begin{equation}\label{diagr-Brnr}
\xymatrix{
\Brnrr(G/H)\ar@{^{(}->}[d]\ar[r]^{i^{\,\nr}}      &\Brnrr(G_1/H)\ar@{^{(}->}[d] \\
 \Brr(G/H)\ar[r]^{i^*}                            &\Brr(G_1/H).
}
\end{equation}

Note that the map $i\colon G_1\to G$ in diagram \eqref{eq:diagr-G-H} is an $H$-equivariant map
from the right $H$-torsor $G_1$ over $G_1/H$ to the right $H$-torsor $G$ over $G/H$.
Sansuc's exact sequence \cite[(6.10.1)]{Sansuc}, applied to this diagram,
gives a commutative diagram with exact rows
\begin{equation}\label{eq:diagr-Sansuc}
\xymatrix{
0=\Pic G \ar[r] &\Pic H\ar[r]\ar[d]^{\id}  &\Brr(G/H)\ar[r]\ar[d]^{i^*}  &\Br G\ar[d] \\
0=\Pic G_1\ar[r] &\Pic H\ar[r]              &\Brr(G_1/H)\ar[r]            &\Br G_1.
}
\end{equation}
We see from  \eqref{eq:diagr-Sansuc}
that the homomorphism $i^*$ restricted to $\ker[\Brr(G/H) \to\Br G]$ is injective,
and we see from \eqref{eq:rational-group} that $i^*$ restricted to $\Brnrr(G/H)$ is injective.
Now it follows from diagram \eqref{diagr-Brnr} that the homomorphism
$$
i^{\,\nr}\colon \Brnrr(G/H)\to \Brnrr(G_1/H)
$$
is injective.
Since $G_1$ is simply connected and $H$ is connected,
by  Bogomolov's theorem \cite[Thm.~9.13]{Mumbai} we have $\Brnrr(G_1/H)=0$.
We conclude that $\Brnrr(G/H)=0$, which proves Theorem \ref{thm:main}.
\qed

\end{document}